\documentclass{birkjour}
\usepackage{amsmath, amsthm, amscd, amsfonts, amssymb, graphicx, color}
\DeclareMathOperator{\Tr}{Tr}
 \newtheorem{thm}{Theorem}[section]
 \newtheorem{cor}[thm]{Corollary}

 \theoremstyle{definition}
 \newtheorem{defn}[thm]{Definition}
 \theoremstyle{remark}
 \newtheorem{rem}[thm]{Remark}
 
 \numberwithin{equation}{section}

\begin{document}

%
%
%
%
%
%
%
%
%

\title[]
{$\ast$-$\eta$-Ricci-Yamabe solitons on  $\alpha$-Cosymplectic manifolds with a quarter-symmetric metric connection }

\author[S. Roy]{Soumendu Roy}

\address{Department of Mathematics,Jadavpur University, Kolkata-700032, India}

\email{soumendu1103mtma@gmail.com}

\thanks{The first author is the corresponding author, supported by Swami Vivekananda Merit Cum Means Scholarship, Government of West Bengal, India.}
\author[S. Dey]{Santu Dey}
\address{Department of Mathematics, Bidhan Chandra College, Asansol - 4, West Bengal-713304 , India}
\email{santu.mathju@gmail.com}
\author[A. Bhattacharyya]{Arindam Bhattacharyya}

\address{Department of Mathematics,Jadavpur University, Kolkata-700032, India}

\email{bhattachar1968@yahoo.co.in}
\author[M. D. Siddiqi]{Mohd. Danish Siddiqi}

\address{Department of Mathematics, Faculty of Science, Jazan University, 82715, Jazan, Kingdom of Saudi Arabia.}

\email{ msiddiqi@jazanu.edu.sa}



\begin{abstract}
The goal of the present paper is to deliberate certain types of metric such as $*$-$\eta$-Ricci-Yamabe soliton on $\alpha$-Cosymplectic manifolds with respect to quarter-symmetric metric connection. Further, we have proved some curvature properties of $\alpha$-Cosymplectic manifolds admitting quarter-symmetric metric connection. Here, we have shown the characteristics of the soliton when the manifold satisfies quarter-symmetric metric connection on $\alpha$-Cosymplectic manifolds. Later, we have acquired Laplace equation from $*$-$\eta$-Ricci-Yamabe soliton equation when the potential vector field $\xi$ of the soliton is of gradient type in terms of quarter-symmetric metric connection. Next, we have developed the nature of the soliton when the vector field is conformal killing admitting quarter-symmetric metric connection. Finally, we present an example of a 5-dimensional $\alpha$-cosymplectic metric as a $*$-$\eta$-Ricci-Yamabe soliton with respect to a quarter-symmetric metric connection to prove our results.
\end{abstract}

\maketitle
\section{Introduction}
The concept of Ricci flow, which is an evolution equation for metrics on a Riemannian manifold, was introduced by R. S. Hamilton \cite{rsham} in the year 1982. A self-similar solution to the Ricci flow \cite{rsham}, \cite{topping} is called a Ricci soliton \cite{rsha} if it moves only by a one parameter family of diffeomorphism and scaling. After the introction of Ricci flow, to construct Yamabe metrics on compact Riemannian manifolds, Hamilton \cite{rsha} developed the idea of Yamabe flow.\\
In 2-dimension the Yamabe flow is equivalent to the Ricci flow \cite{rsham}. But in dimension $> 2$ the Yamabe and Ricci flows do not agree, since the Yamabe flow preserves the conformal class of the metric but the Ricci flow does not in general.\\
Since the introduction of these geometric flows, the respective solitons and their generalizations have been a great centre of attention of many geometers viz. \cite{dong, joma, soumendu,roy,roy2, roy3, roy1, roy4, roy5, roy6,  rshar, sharma, mohd, Dan,singh} who have provided new approaches to understand the geometry of different kinds of Riemannian manifold.\\
Recently in 2019, S. G\"uler and M. Crasmareanu \cite{guler} introduced a new geometric flow which is a scalar combination of Ricci and Yamabe flow under the name Ricci-Yamabe map. This flow is also known as Ricci-Yamabe flow of the type $(\rho, q)$.\\
Let $M$ be a Riemannian manifold of dimension $n$ and $T^{s}_2 (M)$ be the linear space of its symmetric tensor fields of (0, 2)-type and $Riem(M) \subsetneqq T^{s}_2 (M)$  be the infinite space of its Riemannian metrics. In \cite{guler}, the authors have stated the following definition:
\begin{defn}
\cite{guler}
A Riemannian flow on $M$ is a smooth map:
 $$g:I\subseteq \mathbb{R}\rightarrow Riem(M),$$
 where $I$ is a given open interval. We can call it also as time-dependent (or non-stationary) Riemannian metric.
\end{defn}
\begin{defn}
\cite{guler} The map $RY^{(\rho, q, g)}: I \rightarrow T^{s}_2 (M)$ given by:
$$ RY^{(\rho, q, g)}:= \frac{\partial }{\partial t}g(t)+2\rho S(t)+q r(t)g(t),$$
is called the $(\rho, q)$-Ricci-Yamabe map of the Riemannian flow of the Riemannian flow $(M^n, g)$, where $\rho, q$ are some scalars. If $RY^{(\rho, q, g)} \equiv 0$, then $g(\cdot)$ will be called an $(\rho, q)$-Ricci-Yamabe flow.\\
Also in \cite{guler}, the authors characterized that the $(\rho, q)$-Ricci-Yamabe flow is said to be,\\
$\bullet$ Ricci flow \cite{rsham} if $\rho = 1$, $q = 0$.\\
$\bullet$ Yamabe flow \cite{rsha} if $\rho = 0$, $q = 1$.\\
$\bullet$ Einstein flow \cite{CATINO} if $\rho = 1$, $q = -1$.\\\\
A soliton to the Ricci-Yamabe flow is called Ricci-Yamabe solitons\cite{mohd1} if it moves only by one parameter group of diffeomorphism and scaling. The metric of the Riemannain manifold $(M^n, g)$, $n >2$ is said to admit $(\rho, q)$--Ricci-Yamabe soliton or simply Ricci-Yamabe soliton (RYS) $(g, V, \Lambda,\rho, q)$ if it satisfies the equation,
\begin{equation}\label{1.1}
  \pounds_V g+2\rho S=[2\Lambda-q r]g,
\end{equation}
where $\pounds_V g$ denotes the Lie derivative of the metric $g$ along the vector field $V$, $S$ is the Ricci tensor, $r$ is the scalar curvature and $\Lambda, \rho,q$ are real scalars.
\end{defn}
Moreover the Ricci-Yamabe soliton is said to be expanding, steady or shrinking according as $\Lambda$ is negative, zero, positive respectively. Also if $\Lambda, \rho, q$ become smooth function then \eqref{1.1} is called almost Ricci-Yamabe soliton.\\
In 2020, Mohd. Danish Siddiqi et al. \cite{mohd1} introduced a new generalization of Ricci-Yamabe soliton, namely $\eta$-Ricci-Yamabe soliton, which is given by,
\begin{equation}\label{1.2}
  \pounds_V g+2\rho S+[2\Lambda-q r]g+2\mu \eta \otimes \eta=0,
\end{equation}
where $\mu$ is a constant and $\eta$ is a 1-form on $M$.\\
Recently $\ast$-$\eta$ Ricci soliton, a generalization of $\eta$-Ricci soliton \cite{joma}, has been defined by S. Dey and S. Roy \cite{roy6}, which can be given as,
\begin{equation}\label{1.3}
  \pounds_\xi g + 2S^\ast + 2\Lambda g + 2 \mu \eta \otimes \eta =0,
\end{equation}
where $S^\ast$ is the $*$-Ricci tensor, which was first introduced by Tachibana \cite{tachi} on almost Hermitian manifolds and then studied by Hamada \cite{hama} on real hypersurfaces of non-flat complex space forms.\\
Motivating from the above generalizations, we now introduce the notion of $\ast$-$\eta$-Ricci-Yamabe soliton as:
\begin{defn}
A Riemannian manifold $(M,g)$ of dimension $n$ is said to admit $\ast$-$\eta$-Ricci-Yamabe soliton if
\begin{equation}\label{1.4}
\pounds_\xi g+2\rho S^\ast+[2\Lambda-q r^\ast]g+2\mu \eta \otimes \eta=0,
\end{equation}
where $r^\ast=\Tr(S^\ast)$, is the $\ast$-scalar curvature.
\end{defn}
Moreover, in the above equation, if the vector field $\xi$ is the gradient of a smooth function $f$ (denoted by $Df$, $D$ denotes the gradient operator) then the equation \eqref{1.4} is called gradient $*$-$\eta$-Ricci-Yamabe soliton and it is defined as:
\begin{equation}\label{1.5}
  Hess f+\rho S^\ast+\big[\Lambda-\frac{q{r^\ast}}{2}\big] g+\mu \eta \otimes \eta =0,
\end{equation}
where $Hess f$ is the Hessian of the smooth function $f$.\\\\
$~~~~~~~~$ On the other hand, A linear connection $\tilde{\nabla}$ on an $n$-dimensional Riemannian manifold $(M,g)$ is called a quarter-symmetric connection \cite{golab,de1} if its torsion tensor of the connection $\tilde{\nabla}$,
\begin{equation}\label{1.6}
  T(X,Y)=\tilde{\nabla}_X Y-\tilde{\nabla}_Y X-[X,Y]
\end{equation}
satisfies,
\begin{equation}\label{1.7}
  T(X,Y)=\eta(Y)\phi X-\eta(X)\phi Y,
\end{equation}
where $\eta$ is a differentiable 1-form and $\phi$ is a (1, 1) tensor field.\\
If moreover the connection $\tilde{\nabla}$ satisfies,
$$(\tilde{\nabla}_X g)(Y,Z)=0,$$
for ali vector fields $X,Y,Z$ on $(M,g)$, then it is called a quarter-symmetric metric connection \cite{golab,de1}.\\
If $\phi X$ is replaced by $X$, then the connection is called a semi-symmetric metric connection \cite{yano*}.\\
Based on the above facts and discussions in the research of contact geometry, a natuaral \textbf{question} arises \\
{\em Are there contact metric almost manifolds, whose metrics are $*$-$\eta$-Ricci-Yamabe soliton?}\\
In some next sections, we show that indeed the answer to this question is in affirmative. This paper is organized as follows: After introduction, in section 2, we have discussed some basic formulas of contact geometry and curvature properties of $\alpha$-cosymplectic manifolds. In Section 3, we have proved some characteristics of $n$-dimensional $\alpha$-cosymplectic manifold admitting quarter-symmetric metric connection. Section 4 deals with $*$-$\eta$-Ricci-Yamabe soliton in terms of quarter-symmetric metric connection on $\alpha$-cosymplectic manifolds. In this section, we have established the relationship between soliton constans $\Lambda$ and $\mu$. Section 4 is also devoted to form the Laplace equation from $\eta$-Ricci-Yamabe soliton equation with respect to quarter-symmetric metric connection when the potential vector field $\xi$ of the soliton is of gradient type.  We have also constructed some applications on Laplace equation and proved some theorems on harmonic aspect of $\ast$-$\eta$-Ricci-Yamabe soliton on $n$-dimensional $\alpha$-cosymplectic manifold with a quarter-symmetric metric connection in section 5 and 6 respectively. We have next considered the potential vector field $V$ of the solition as conformal Killing vector field to characterize the vector field to accessorize the nature of this soliton on this manifold with quarter-symmetric metric connection. In section 7, we have constructed an example to illustrate the existence of $*$-$\eta$-Ricci-Yamabe soliton on $5$-dimensional $\alpha$-cosymplectic manifold with restect to quarter-symmetric metric connection.

\section{Preliminaries}
Let $M$ be a (2$n$+1)dimensional connected almost contact metric manifold with an almost contact metric structure $(\phi, \xi, \eta, g)$ where $\phi$ is a $(1,1)$ tensor field, $\xi$ is a vector field, $\eta$ is a 1-form  and $g$ is the compatible Riemannian metric such that \cite{blair},
\begin{equation}\label{2.1}
\phi^2(X) = -X + \eta(X)\xi, \eta(\xi) = 1, \eta \circ \phi = 0, \phi \xi = 0,
\end{equation}\\
\begin{equation}\label{2.2}
g(\phi X,\phi Y) = g(X,Y) - \eta(X)\eta(Y),
\end{equation}\\
\begin{equation}\label{2.3}
g(X,\phi Y) = -g(\phi X,Y),
\end{equation}\\
\begin{equation}\label{2.4}
g(X,\xi) = \eta(X),
\end{equation}\\
for all vector fields $X, Y \in \chi(M),$ where $\chi(M)$ denotes the collection of all smooth vector fields of $M$.\\
On this manifold $M$, the 2-form $\Phi$ is defined as,
\begin{equation}\label{2.5}
  \Phi(X,Y)=g(\phi X,Y),
\end{equation}
for all vector fields $X, Y \in \chi(M).$\\
An almost contact metric manifold $(M, \phi, \xi, \eta, g)$ is said to be almost cosymplectic \cite{goldbarg,has} if $d\eta = 0$ and $d\Phi = 0$, where $d$ is the exterior differential operator. An almost contact manifold $(M, \phi, \xi, \eta, g)$ is said to be normal if the Nijenhuis
torsion,
$$N_\phi (X,Y)=[\phi X,\phi Y-\phi[\phi X,Y]-\phi[X,\phi Y]+\phi^2[X,Y] + 2d\eta(X,Y)\xi,$$
vanishes for any vector fields $X$ and $Y$ . A normal almost cosymplectic manifold is called a cosymplectic manifold \cite{goldbarg,has}.\\
An almost contact metric manifold $M$ is said to be almost $\alpha$-Kenmotsu if $d\eta = 0$ and $d\Phi = 2\alpha\eta\wedge\Phi$, $\alpha$ being a non-zero real constant.\\
In 2005, a combination of almost $\alpha$-Kenmotsu and almost cosymplectic manifolds was developed by Kim and Pak \cite{kim}, into a new class called almost $\alpha$-cosymplectic manifold, where $\alpha$ is a scalar. An almost $\alpha$-cosymplectic manifold can be defined as \cite{hakan, has},
$$d\eta=0, \quad d\Phi = 2\alpha\eta\wedge\Phi,$$
for any real number $\alpha$.\\
A normal almost $\alpha$-cosymplectic manifold is called an $\alpha$-cosymplectic manifold. An $\alpha$-cosymplectic manifold is either cosymplectic under the condition $\alpha = 0$ or $\alpha$-Kenmotsu under the condition $\alpha \neq 0$, for $\alpha \in \mathbb{R}$ \cite{has}. \\
In an $\alpha$-cosymplectic manifold, we have \cite{hakan, has},
\begin{equation}\label{2.6}
  (\nabla_X \phi)Y=\alpha(g(\phi X, Y)\xi, \eta(Y)\phi X).
\end{equation}
Let $M$ be an $n$-dimensional $\alpha$-cosymplectic manifold. Then from \eqref{2.6}, we have,
\begin{equation}\label{2.7}
  \nabla_X \xi = -\alpha \phi^2X = \alpha[X - \eta(X)\xi],
\end{equation}
where $\nabla$ is the Levi-Civita connection associated with $g$.\\
On an $n$-dimensional $\alpha$-cosymplectic manifold $M$, the following relations hold \cite{has},
\begin{equation}\label{2.8}
  R(X,Y)\xi=\alpha^2[\eta(X)Y-\eta(Y)X],
\end{equation}
\begin{equation}\label{2.9}
   R(\xi, X)Y=\alpha^2[\eta(Y)X-g(X,Y)\xi],
\end{equation}
\begin{equation}\label{2.10}
  R(\xi, X)\xi=\alpha^2[X-\eta(X)\xi],
\end{equation}e
\begin{equation}\label{2.11}
  \eta(R(X,Y)Z)=\alpha^2[\eta(Y)g(X,Z)-\eta(X)g(Y,Z)],
\end{equation}
\begin{equation}\label{2.12}
  S(X,\xi)=-\alpha^2(n-1)\eta(X),
\end{equation}
for all vector fields $X, Y, Z$ on $M$, where $R$ is the Riemannian curvature tensor and $S$ is the Ricci tensor of $M$.\\
In paper \cite{has}, Lemma 2.2, the authors have proved that the $\ast$-Ricci tensor on an $n$-dimensional $\alpha$-cosymplectic manifold is given by,
\begin{equation}\label{2.13}
  S^\ast(Y,Z)=S(Y,Z)+\alpha^2(n-2)g(Y,Z)+\alpha^2\eta(Y)\eta(Z),
\end{equation}
for any vector field $Y,Z$ on $M$, where $S$ and $S^\ast$ are the Ricci tensor and the $\ast$-Ricci tensor of type $(0, 2)$, respectively
on $M$.\\
Taking $Y=Z=e_i$, where $e_i's$ are the orthonormal basis of $T_p(M)$ for $i=1,2,...n$, we get,
\begin{equation}\label{2.14}
  r^\ast=r+\alpha^2(n-1)^2,
\end{equation}
where $r^\ast=\Tr(S^\ast)$, is the $\ast$-scalar curvature and $r$ is the scalar curvature.\\
\section{On an $n$-dimensional $\alpha$-cosymplectic manifold with a quarter-symmetric metric connection}
Let $\tilde{\nabla}$ be a linear connection and $\nabla$ be a Levi-Civita connection of an almost
contact metric manifold $M$ such that,
\begin{equation}\label{3.1}
  \tilde{\nabla}_X Y=\nabla_X Y+H(X,Y),
\end{equation}
where $H$ is a tensor of type $(1, 1)$. For $\tilde{\nabla}$ to be a quarter-symmetric metric connection on $M$, we have \cite{golab},
\begin{equation}\label{3.2}
  H(X,Y)=\frac{1}{2}\big[T(X,Y)+T^\prime(X,Y)+T^\prime(Y,X)\big],
\end{equation}
where
\begin{equation}\label{3.3}
  g(T^\prime(X,Y),Z)=g(T(Z,X),Y).
\end{equation}
Now from \eqref{1.7} and \eqref{3.3}, we get,
\begin{equation}\label{3.4}
  T^\prime(X,Y)=g(\phi Y,X)\xi-\eta(X)\phi Y.
\end{equation}
Using \eqref{1.6}, \eqref{3.4} and \eqref{3.2}, we obtain,
\begin{equation}\label{3.5}
  H(X,Y)=-\eta(X)\phi Y.
\end{equation}
Hence a quarter-symmetric metric connection $\tilde{\nabla}$ in a $\alpha$-cosymplectic manifold is given by,
\begin{equation}\label{3.6}
  \tilde{\nabla}_X Y=\nabla_X Y-\eta(X)\phi Y.
\end{equation}
Let $R$ and $\tilde{R}$ be the curvature tensors of $\nabla$ and $\tilde{\nabla}$ of a $\alpha$-cosymplectic manifold, respectively.\\
Then using \eqref{3.6} and \eqref{2.7}, we have,
\begin{equation}\label{3.7}
  \tilde{R}(X,Y)Z=R(X,Y)Z+\eta(X)(\nabla_Y \phi)Z-\eta(Y)(\nabla_X \phi)Z,
\end{equation}
which in view of \eqref{2.6}, gives,
\begin{eqnarray}\label{3.8}
  \tilde{R}(X,Y)Z&=&R(X,Y)Z+\alpha\eta(X)g(\phi Y, Z)\xi-\alpha\eta(Y)g(\phi X, Z)\xi \nonumber\\
                 &-&\alpha\eta(X)\eta(Z)\phi Y+\alpha\eta(Y)\eta(Z)\phi X.
\end{eqnarray}
A relation between the curvature tensor of M with respect to the quarter-symmetric metric connection $\tilde{\nabla}$ and Levi-Civita connection $\nabla$ is given by the equation \eqref{3.8}.\\
Taking inner product of \eqref{3.8} with $W$, we get,
\begin{eqnarray}\label{3.9}
  \tilde{R}(X,Y,Z,W) &=& R(X,Y,Z,W)+\alpha\eta(X)\eta(W)g(\phi Y,Z) \nonumber \\
                     &-& \alpha\eta(Y)\eta(W)g(\phi X, Z)-\alpha\eta(X)\eta(Z)g(\phi Y,W) \nonumber \\
                     &+& \alpha\eta(Y)\eta(Z)g(\phi X,W).
\end{eqnarray}
Contracting \eqref{3.9} over $X$ and $W$, we obtain,
\begin{equation}\label{3.10}
  \tilde{S}(Y,Z)=S(Y,Z)+\alpha g(\phi Y,Z),
\end{equation}
where $\tilde{S}$ and $S$ are the Ricci tensors of the connections $\tilde{\nabla}$ and $\nabla$, respectively.\\
So in a $\alpha$-cosymplectic manifold, the Ricci tensor of the quarter-symmetric metric connection is not symmetric.\\
Again, contracting \eqref{3.10} over $Y,Z$, it gives,
\begin{equation}\label{3.11}
  \tilde{r}=r.
\end{equation}
where $\tilde{r}$ and $r$ are the scalar curvatures of the connections $\tilde{\nabla}$ and $\nabla$, respectively.\\
So we have the following theorem:
\begin{thm}
For a $\alpha$-cosymplectic manifold $M$ with the quarter-symmetric metric connection $\tilde{\nabla}$\\
(i) The Riemannian curvature tensor $\tilde{R}$ is given by \eqref{3.8},\\
(ii) The Ricci tensor $\tilde{S}$ is given by \eqref{3.10},\\
(iii) $\tilde{R}(X, Y,Z, W) + \tilde{R}(X, Y, W, Z) = 0$,\\
(iv) $\tilde{R}(X, Y, Z, W) + \tilde{R}(Y, X, Z, W) = 0$,\\
(v) $\tilde{S}(Y, \xi) = S(Y, \xi) = -\alpha^2(n-1)\eta(Y)$,\\
(vi) The Ricci tensor $\tilde{S}$ is not symmetric,\\
(vii)  $\tilde{r}= r$.
\end{thm}
\section{Main Results}
Let us consider an $n$-dimensional $\alpha$-cosymplectic manifold $M$, admitting a $\ast$-$\eta$-Ricci-Yamabe soliton.\\
Then from \eqref{1.4}, we get,
\begin{equation}\label{4.1}
  (\pounds_\xi g)(Y,Z)+2\rho S^\ast(Y,Z)+[2\Lambda-q r^\ast]g(Y,Z)+2\mu \eta(Y)\eta(Z)=0,
\end{equation}
for any vector fields $Y,Z$ on $M$.\\
Using \eqref{2.13} and \eqref{2.14}, the above equation takes the form,
\begin{multline}\label{4.2}
  (\pounds_\xi g)(Y,Z)+2\rho S(Y,Z)+\big[2\Lambda+2\alpha^2\rho(n-2)-q r-q\alpha^2(n-1)^2\big]g(Y,Z)\\
  +\big[2\mu+2\alpha^2\rho\big] \eta(Y)\eta(Z)=0.
\end{multline}
We consider $\tilde{\pounds}_V$ is the Lie derivative along $V$ with respect to the quarter-symmetric metric connection $\tilde{\nabla}$.\\
Now if an $n$-dimensional $\alpha$-cosymplectic manifold $M$ admits a $\ast$-$\eta$-Ricci-Yamabe soliton with respect to the quarter-symmetric metric connection $\tilde{\nabla}$, then from \eqref{4.2}, we have,
\begin{multline}\label{4.3}
  (\tilde{\pounds}_\xi g)(Y,Z)+2\rho \tilde{S}(Y,Z)+\big[2\Lambda+2\alpha^2\rho(n-2)-q \tilde{r}-q\alpha^2(n-1)^2\big]g(Y,Z)\\
  +\big[2\mu+2\alpha^2\rho\big] \eta(Y)\eta(Z)=0.
\end{multline}
In view of \eqref{3.6}, we obtain,
\begin{eqnarray}\label{4.4}
  (\tilde{\pounds}_V g)(X,Y) &=& g(\tilde{\nabla}_X V,Y)+g(X,\tilde{\nabla}_Y V) \nonumber \\
                             &=& (\pounds_V g)(X,Y)-\eta(X)g(\phi V,Y)-\eta(Y)g(X, \phi V).\nonumber\\
 \end{eqnarray}
Then using \eqref{2.1}, we get,
\begin{equation}\label{4.5}
  (\tilde{\pounds}_\xi g)(X,Y)=(\pounds_\xi g)(X,Y),
\end{equation}
for any vector fields $X,Y$ on $M$.\\
From \eqref{4.5}, \eqref{3.10} and \eqref{3.11}, \eqref{4.3} becomes,
\begin{multline}\label{4.6}
  \Big[(\pounds_\xi g)(Y,Z)+2\rho S(Y,Z)+\big[2\Lambda+2\alpha^2\rho(n-2)-q r-q\alpha^2(n-1)^2\big]g(Y,Z)\\
  +\big[2\mu+2\alpha^2\rho\big] \eta(Y)\eta(Z)\Big]+2\alpha\rho g(\phi Y,Z)=0. \\
\end{multline}
Now in an $n$-dimensional $\alpha$-cosymplectic manifold, from \eqref{2.7}, we obtain,
\begin{equation}\label{4.7}
(\pounds_\xi g)(Y,Z)=g(\nabla_Y \xi,Z)+g(Y, \nabla_Z \xi)=2\alpha[g(Y,Z)-\eta(Y)\eta(Z)].
\end{equation}
Putting this value in \eqref{4.6} and then taking $Z=\xi$, we get,
\begin{equation}\label{4.8}
  2\rho S(Y,\xi)+\big[2\Lambda+2\alpha^2\rho(n-2)-q r-q\alpha^2(n-1)^2\big]\eta(Y)+\big[2\mu+2\alpha^2\rho\big] \eta(Y)=0.
\end{equation}
Using \eqref{2.12}, we find,
\begin{equation}\label{4.9}
  \Lambda+\mu=\frac{qr}{2}+\frac{q\alpha^2(n-1)^2}{2}.
\end{equation}
Hence we can state,
\begin{thm}
If the metric of an $n$-dimensional $\alpha$-cosymplectic manifold admits a $\ast$-$\eta$-Ricci-Yamabe soliton with respect to a quarter-symmetric metric connection $\tilde{\nabla}$, then the soliton constants $\Lambda$ and $\mu$ are related by the equation \eqref{4.9}.
\end{thm}
Again using \eqref{4.2}, \eqref{4.6} takes the form,
\begin{equation}\label{4.10}
  \alpha\rho g(\phi Y,Z)=0.
\end{equation}
for any vector fields $Y,Z$ on $M$.\\
Then $\alpha=0$, since in a $\ast$-$\eta$-Ricci-Yamabe soliton $\rho$ can not be zero.\\
So as $\alpha=0$, then $M$ reduces to an $n$-dimensional cosymplectic manifold.\\
This gives,
\begin{cor}
If the metric of an $n$-dimensional $\alpha$-cosymplectic manifold admits a $\ast$-$\eta$-Ricci-Yamabe soliton with respect to the Levi-Civita connection $\nabla$. Then the metric admits the soliton with respect to a quarter-symmetric metric connection $\tilde{\nabla}$ iff the manifold becomes an $n$-dimensional cosymplectic manifold.
\end{cor}
Taking $Y=Z=e_i$, where $e_i's$ are the orthonormal basis of $T_p(M)$ for $i=1,2,...n$, we get,
\begin{equation}\label{4.11}
  div \xi+\rho r+n\Big[\Lambda+\alpha^2\rho(n-2)-\frac{q r}{2} -\frac{q\alpha^2(n-1)^2}{2} \Big]+\mu+\alpha^2\rho=0,
\end{equation}
where $div \xi$ is the divergence of the vector field $\xi$.\\
Replacing the value of $\mu$ from \eqref{4.9} in the above equation, $\Lambda$ takes the form,
\begin{equation}\label{4.12}
  \Lambda=\frac{qr}{2}-\frac{q\alpha^2(n-1)^2}{2}-\alpha^2\rho(n-1)-\frac{div \xi+\rho r}{n-1}.
\end{equation}
And in view of \eqref{4.12}, \eqref{4.9} becomes,
\begin{equation}\label{4.13}
  \mu=\alpha^2\rho(n-1)+\frac{div \xi+\rho r}{n-1}.
\end{equation}
Hence we can state,
\begin{cor}
If the metric of an $n$-dimensional $\alpha$-cosymplectic manifold admits a $\ast$-$\eta$-Ricci-Yamabe soliton with respect to a quarter-symmetric metric connection $\tilde{\nabla}$. Then the soliton constants $\Lambda$ and $\mu$ takes the form of $\frac{qr}{2}-\frac{q\alpha^2(n-1)^2}{2}-\alpha^2\rho(n-1)-\frac{div \xi+\rho r}{n-1}$ and $\alpha^2\rho(n-1)+\frac{div \xi+\rho r}{n-1}$ respectively, where $div \xi$ is the divergence of the vector field $\xi$.
\end{cor}
If $\xi=grad(f)$, where $grad(f)$ is the gradient of a smooth function $f$, then from \eqref{4.11}, we have the following,
\begin{thm}\label{TD1}
Let an $n$-dimensional $\alpha$-cosymplectic manifold admit a $\ast$-$\eta$-Ricci-Yamabe soliton with respect to a quarter-symmetric metric connection $\tilde{\nabla}$. If the vector field $\xi$, associated with the soliton is of the form $grad(f)$, where $f$ is a smooth function. Then the Laplacian equation satisfied by $f$ becomes:
\begin{equation}\label{4.14}
  \Delta(f)=-n\Big[\Lambda+\alpha^2\rho(n-2)-\frac{q r}{2} -\frac{q\alpha^2(n-1)^2}{2} \Big]-\mu-\rho(\alpha^2+r).
\end{equation}
\end{thm}
\section{Some Applications}
As an application, we obtain the following results for $\ast$-$\eta$-Ricci soliton, $\ast$-$\eta$-Yamabe soliton, and $\ast$-$\eta$-Einstein soliton ($\rho=1, q=0$, $\rho=0, q=1$, and $\rho=1, q=-1$ ( cf. \cite{rsham, rsha, CATINO})).

\begin{thm}\label{TD2}
Let an $n$-dimensional $\alpha$-cosymplectic manifold admit a $\ast$-$\eta$-Ricci soliton with respect to a quarter-symmetric metric connection $\tilde{\nabla}$. If the vector field $\xi$, associated with the soliton is of the form $grad(f)$, where $f$ is a smooth function. Then the Laplacian equation satisfied by $f$ becomes:
\begin{equation}\label{d1}
  \Delta(f)=-n\Big[\Lambda+\alpha^2(n-2) \Big]-\mu-(\alpha^2+r).
\end{equation}
\end{thm}
\begin{thm}\label{TD3}
Let an $n$-dimensional $\alpha$-cosymplectic manifold admit a $\ast$-$\eta$-Yamabe soliton with respect to a quarter-symmetric metric connection $\tilde{\nabla}$. If the vector field $\xi$, associated with the soliton is of the form $grad(f)$, where $f$ is a smooth function. Then the Laplacian equation satisfied by $f$ becomes:
\begin{equation}\label{d2}
  \Delta(f)=-n\Big[\Lambda-\frac{ r}{2} -\frac{\alpha^2(n-1)^2}{2} \Big]-\mu.
\end{equation}
\end{thm}
\begin{thm}\label{TD4}
Let an $n$-dimensional $\alpha$-cosymplectic manifold admit a $\ast$-$\eta$-Einstein soliton with respect to a quarter-symmetric metric connection $\tilde{\nabla}$. If the vector field $\xi$, associated with the soliton is of the form $grad(f)$, where $f$ is a smooth function. Then the Laplacian equation satisfied by $f$ becomes:
\begin{equation}\label{d3}
  \Delta(f)=-n\Big[\Lambda+\alpha^2(n-2)+\frac{ r}{2} +\frac{\alpha^2(n-1)^2}{2} \Big]-\mu-(\alpha^2+r).
\end{equation}
\end{thm}
\section{Harmonic aspect of $\ast$-$\eta$-Ricci-Yamabe soliton on $n$-dimensional $\alpha$-cosymplectic manifold with a quarter-symmetric metric connection}

This section is based on  the fact that a function $f:{M}\longrightarrow \mathbb{R}$ is said to be harmonic if $\Delta f=0$, where $\Delta$ is the Lalplacian operator on $\textbf{M}$ \cite{Yau}. Therefore, let us considering the fact that the vector field $\xi$ is a gradient of a harmonic function  $f$, then  from Theorem (\ref{TD1}), we turn up the following conclusions:\\

\begin{thm}\label{TH1}
Let an $n$-dimensional $\alpha$-cosymplectic manifold $M$ admit a $\ast$-$\eta$-Ricci-Yamabe soliton with respect to a quarter-symmetric metric connection $\tilde{\nabla}$ and the vector field $\xi$, associated with the soliton is of the form $grad(f)$,
 if $f$ is a harmonic function  on  ${M}$.  Then the $\ast$-$\eta$-Ricci-Yamabe soliton is expanding, steady and shrinking according as
	\begin{enumerate}
\item \quad $\frac{q}{2}\left\{r+\alpha^{2}(n-1)^{2}\right\} > \left\{\mu+\rho(\alpha^{2}+r)+\alpha^{2}\rho(n-2)\right\}$,\\
\item\quad $\frac{q}{2}\left\{r+\alpha^{2}(n-1)^{2}\right\} = \left\{\mu+\rho(\alpha^{2}+r)+\alpha^{2}\rho(n-2)\right\}$ \quad and \\
\item\quad $\frac{q}{2}\left\{r+\alpha^{2}(n-1)^{2}\right\} < \left\{\mu+\rho(\alpha^{2}+r)+\alpha^{2}\rho(n-2)\right\}$\quad respectively.
\end{enumerate}
	\end{thm}
	\begin{proof}
	Form equation (\ref{4.14}) we can easily obtain the desired result.
	\end{proof}
\begin{thm}\label{TH2}
Let an $n$-dimensional $\alpha$-cosymplectic manifold $M$ admit a $\ast$-$\eta$-Ricci soliton with respect to a quarter-symmetric metric connection $\tilde{\nabla}$ and the vector field $\xi$, associated with the soliton is of the form $grad(f)$,
 if $f$ is a harmonic function  on  ${M}$.  Then the $\ast$-$\eta$-Ricci soliton is shrinking.
	\end{thm}
	\begin{thm}\label{TH3}
Let an $n$-dimensional $\alpha$-cosymplectic manifold $M$ admit a $\ast$-$\eta$-Yamabe soliton with respect to a quarter-symmetric metric connection $\tilde{\nabla}$ and the vector field $\xi$, associated with the soliton is of the form $grad(f)$,
 if $f$ is a harmonic function  on  ${M}$.  Then the $\ast$-$\eta$-Yamabe soliton is expanding, steady and shrinking according as
	\begin{enumerate}
\item \quad $\frac{1}{2}\left\{r+\alpha^{2}(n-1)^{2}\right\} > \mu$,\\
\item\quad $\frac{1}{2}\left\{r+\alpha^{2}(n-1)^{2}\right\} = \mu $ \quad and \\
\item\quad $\frac{1}{2}\left\{r+\alpha^{2}(n-1)^{2}\right\} < \mu $\quad respectively.
\end{enumerate}
	\end{thm}
	\begin{proof}
	Form equation (\ref{4.14}) we can easily obtain the desired result.
	\end{proof}
	\begin{thm}\label{TH1}
Let an $n$-dimensional $\alpha$-cosymplectic manifold $M$ admit a $\ast$-$\eta$-Einstein soliton with respect to a quarter-symmetric metric connection $\tilde{\nabla}$ and the vector field $\xi$, associated with the soliton is of the form $grad(f)$,
 if $f$ is a harmonic function  on  ${M}$.  Then the $\ast$-$\eta$-Einstein soliton is also shrinking.
	\end{thm}
	\begin{proof}
	Form equation (\ref{4.14}) we can easily obtain the desired result.
	\end{proof}
\begin{defn}
A vector field $V$ is said to be a conformal Killing vector field iff the following relation holds:
\begin{equation}\label{4.15}
  (\pounds_V g)(Y,Z)=2\theta g(Y,Z),
\end{equation}
where $\theta$ is some function of the co-ordinates(conformal scalar).\\
Moreover if $\theta$ is not constant the conformal Killing vector field $V$ is said to be proper. Also when $\theta$ is constant, $V$ is called homothetic vector field and when the constant $\theta$ becomes non zero, $V$ is said to be proper homothetic vector field. If $\theta = 0$ in the above equation, then $V$ is called Killing vector field.
\end{defn}
Let us consider an $n$-dimensional $\alpha$-cosymplectic manifold $M$, admitting a $\ast$-$\eta$-Ricci-Yamabe soliton, where $\xi=V$.\\
Then using \eqref{1.4}, \eqref{2.13} and \eqref{2.14}, we have,
\begin{multline}\label{4.16}
  (\pounds_V g)(Y,Z)+2\rho S(Y,Z)+\big[2\Lambda+2\alpha^2\rho(n-2)-q r-q\alpha^2(n-1)^2\big]g(Y,Z)\\
  +\big[2\mu+2\alpha^2\rho\big] \eta(Y)\eta(Z)=0.
\end{multline}
Now if an $n$-dimensional $\alpha$-cosymplectic manifold $M$ admits a $\ast$-$\eta$-Ricci-Yamabe soliton $(g,V,\Lambda,\mu,\rho,q)$ with respect to the quarter-symmetric metric connection $\tilde{\nabla}$, then from \eqref{4.16}, we have,
\begin{multline}\label{4.17}
  (\tilde{\pounds}_V g)(Y,Z)+2\rho \tilde{S}(Y,Z)+\big[2\Lambda+2\alpha^2\rho(n-2)-q \tilde{r}-q\alpha^2(n-1)^2\big]g(Y,Z)\\
  +\big[2\mu+2\alpha^2\rho\big] \eta(Y)\eta(Z)=0.
\end{multline}
In the view of \eqref{4.4}, \eqref{3.10}, \eqref{3.11}, the above equation becomes,
\begin{multline}\label{4.18}
  (\pounds_V g)(Y,Z)-\eta(Y)g(\phi V,Z)-\eta(Z)g(Y, \phi V)+2\rho [S(Y,Z)+\alpha g(\phi Y,Z)]\\
  +\big[2\Lambda+2\alpha^2\rho(n-2)-q r-q\alpha^2(n-1)^2\big]g(Y,Z)+\big[2\mu+2\alpha^2\rho\big] \eta(Y)\eta(Z)=0.
\end{multline}
Using \eqref{4.15}, the equation \eqref{4.18} takes the form,
\begin{multline}\label{4.19}
\big[2\theta+2\Lambda+2\alpha^2\rho(n-2)-q r-q\alpha^2(n-1)^2\big]g(Y,Z)+\big[2\mu+2\alpha^2\rho\big] \eta(Y)\eta(Z)\\
  +2\rho [S(Y,Z)+\alpha g(\phi Y,Z)]-\eta(Y)g(\phi V,Z)-\eta(Z)g(Y, \phi V)=0.
\end{multline}
Now, taking $Z=\xi$ and using \eqref{2.12}, we get,
\begin{equation}\label{4.20}
  g(Y,\phi V)=\big[2\theta+2\Lambda+2\mu-qr-q\alpha^2(n-1)^2\big]\eta(Y),
\end{equation}
which can be written as,
\begin{equation}\label{4.21}
 \phi V=\big[2\theta+2\Lambda+2\mu-qr-q\alpha^2(n-1)^2\big]\xi.
\end{equation}
\begin{thm}\label{TCV}
Let an $n$-dimensional $\alpha$-cosymplectic manifold admit a $\ast$-$\eta$-Ricci-Yamabe soliton $(g,V,\Lambda,\mu,\rho,q)$ with respect to a quarter-symmetric metric connection $\tilde{\nabla}$. If $V$ is a conformal Killing vector field, then \eqref{4.20} holds.
\end{thm}

\begin{rem} For particular values of $\rho$ and $q$, we can also obtain the similar results like Theorem (\ref{TCV}) for $\ast$-$\eta$-Ricci soliton, $\ast$-$\eta$-Yamabe soliton, and $\ast$-$\eta$-Einstein soliton, respectively.
\end{rem}
\section{Example of a 5-dimensional $\alpha$-cosymplectic metric as a $\ast$-$\eta$-Ricci-Yamabe soliton with respect to a quarter-symmetric metric connection}
We consider the 5-dimensional manifold $M = \{(x, y, z,u,v) \in \mathbb{R}^5 \}$,
where $(x, y, z, u, v)$ are standard coordinates in $\mathbb{R}^5$. Let ${e_1, e_2, e_3, e_4, e_5}$ be a linearly independent
frame field on M given by,
\begin{equation}
   e_1=e^{\alpha v}\frac{\partial}{\partial x},\quad e_2 = e^{\alpha v}\frac{\partial}{\partial y}, \quad e_3 =e^{\alpha v}\frac{\partial}{\partial z}, \quad e_4=e^{\alpha v}\frac{\partial}{\partial u}, \quad e_5=-\frac{\partial}{\partial v}. \nonumber
\end{equation}
Let $g$ be the Riemannian metric defined by,
\begin{equation}
  g(e_i,e_j)=0,i\neq j,i,j=1,2,3,4,5,\nonumber
\end{equation}
\begin{equation}
   g(e_1,e_1) = g(e_2,e_2) = g(e_3,e_3) =g(e_4,e_4) = g(e_5,e_5)= 1.\nonumber
\end{equation}
Let $\eta$ be the 1-form defined by $\eta(Z) = g(Z,e_5)$, for any $Z \in \chi(M)$,where $\chi(M)$ is the set of all differentiable vector fields on $M$ and $\phi$ be the (1, 1)-tensor field defined by,
\begin{equation}
  \phi e_1=-e_2, \quad \phi e_2=e_1,\quad  \phi e_3=-e_4, \quad \phi e_4= e_3, \quad \phi e_5=0.\nonumber
\end{equation}
Then, using the linearity of $\phi$ and $g$, we have $\eta(e_5) = 1, \phi ^2 (Z) = -Z + \eta(Z)e_5$ and $g(\phi Z,\phi U) = g(Z,U) - \eta(Z)\eta(U)$, for any $Z,U \in \chi(M)$. Thus $e_5 = \xi$.\\
Let $\nabla$ be the Levi-Civita connection with respect to the Riemannian metric $g$. Then we have,\\
  $ [e_1,e_5] =\alpha e_1, [e_2,e_5] =\alpha e_2, [e_3,e_5] =\alpha e_3, [e_4,e_5] =\alpha e_4$, and $ [e_i,e_j] =0$ for others $i,j$.\\
The connection $\nabla$ of the metric $g$ is given by,
\begin{eqnarray}
  2g(\nabla_X Y,Z) &=& Xg(Y,Z)+Yg(Z,X)-Zg(X,Y)\nonumber \\
                   &-& g(X, [Y,Z])-g(Y, [X, Z]) + g(Z, [X, Y ]),\nonumber
\end{eqnarray}
which is known as Koszul’s formula.\\
Using Koszul’s formula, we can easily calculate,
$$\nabla_{e_1} e_1 =-\alpha e_5, \nabla_{e_1} e_2 =0, \nabla_{e_1} e_3 =0 , \nabla_{e_1} e_4=0 ,\nabla_{e_1} e_5=\alpha e_1 ,$$
$$\nabla_{e_2} e_1 =0, \nabla_{e_2} e_2 =-\alpha e_5, \nabla_{e_2} e_3 =0,\nabla_{e_2} e_4=0,\nabla_{e_2} e_5=\alpha e_2,$$
$$\nabla_{e_3} e_1 =0, \nabla_{e_3} e_2 = 0, \nabla_{e_3} e_3 =-\alpha e_5,\nabla_{e_3} e_4=0,\nabla_{e_3} e_5=\alpha e_3,$$
$$\nabla_{e_4} e_1=0,\nabla_{e_4} e_2=0,\nabla_{e_4} e_3=0,\nabla_{e_4} e_4=-\alpha e_5,\nabla_{e_4} e_5=\alpha e_4,$$
$$\nabla_{e_5} e_1=\nabla_{e_5} e_2=\nabla_{e_5} e_3=\nabla_{e_5} e_4=\nabla_{e_5} e_5=0.$$
It can be easily verified that for $e_5=\xi$, the manifold satisfies,\\
$\nabla_X \xi=\alpha[X-\eta(X)\xi]$ and $(\nabla_X \phi)Y=\alpha[g(\phi X,Y)\xi-\eta(Y)\phi X]$.\\
Thus the manifold $M$ is an $\alpha$-cosymplectic manifold.\\
Also, the Riemannian curvature tensor $R$ is given by,
$$R(X, Y )Z = \nabla_X\nabla_Y Z - \nabla_Y \nabla_X Z - \nabla_{[X,Y]} Z.$$
Hence,
$$R(e_1,e_2)e_2 = R(e_1,e_3)e_3 = R(e_1,e_4)e_4 = R(e_1,e_5)e_5=-\alpha^2 e_1,$$
$$R(e_1,e_2)e_1=\alpha^2 e_2, R(e_1,e_3)e_1 = R(e_2,e_3)e_2 = R(e_5,e_3)e_5=\alpha^2 e_3$$
$$R(e_2,e_3)e_3= R(e_2,e_4)e_4= R(e_2,e_5)e_5=-\alpha^2 e_2, R(e_3,e_4)e_4=-\alpha^2 e_3,$$
$$R((e_1,e_5)e_1= R(e_2,e_5)e_2= R(e_3,e_5)e_3= R(e_4,e_5)e_4=\alpha^2 e_5,$$
$$R(e_1,e_4)e_1= R(e_2,e_4)e_2= R(e_3,e_4)e_3= R(e_5,e_4)e_5= \alpha^2 e_4.$$
From the above relations, we have,
$$S(e_1,e_1)= S(e_2,e_2)= S(e_3,e_3)= S(e_4,e_4)= S(e_5,e_5)=-4\alpha^2.$$
Hence $r=-20\alpha^2$.\\
Let $\tilde{\nabla}$ be a quarter-symmetric metric connection on $M$.\\
Then using \eqref{3.6}, we get,
$$\tilde{\nabla}_{e_1} e_1 =-\alpha e_5, \tilde{\nabla}_{e_1} e_2 =0, \tilde{\nabla}_{e_1} e_3 =0 , \tilde{\nabla}_{e_1} e_4=0 ,\tilde{\nabla}_{e_1} e_5=\alpha e_1 ,$$
$$\tilde{\nabla}_{e_2} e_1 =0, \tilde{\nabla}_{e_2} e_2 =-\alpha e_5, \tilde{\nabla}_{e_2} e_3 =0,\tilde{\nabla}_{e_2} e_4=0,\tilde{\nabla}_{e_2} e_5=\alpha e_2,$$
$$\tilde{\nabla}_{e_3} e_1 =0, \tilde{\nabla}_{e_3} e_2 = 0, \tilde{\nabla}_{e_3} e_3 =-\alpha e_5,\tilde{\nabla}_{e_3} e_4=0,\tilde{\nabla}_{e_3} e_5=\alpha e_3,$$
$$\tilde{\nabla}_{e_4} e_1=0,\tilde{\nabla}_{e_4} e_2=0,\tilde{\nabla}_{e_4} e_3=0,\tilde{\nabla}_{e_4} e_4=-\alpha e_5,\tilde{\nabla}_{e_4} e_5=\alpha e_4,$$
$$\tilde{\nabla}_{e_5} e_1=e_2, \tilde{\nabla}_{e_5} e_2=-e_1, \tilde{\nabla}_{e_5} e_3=e_4, \tilde{\nabla}_{e_5} e_4=-e_3, \tilde{\nabla}_{e_5} e_5=0.$$
Then,
$$\tilde{R}(e_1,e_2)e_2 = \tilde{R}(e_1,e_3)e_3 = \tilde{R}(e_1,e_4)e_4 =-\alpha^2 e_1,\tilde{R}(e_1,e_5)e_5= -\alpha^2 e_1-\alpha e_2,$$
$$\tilde{R}(e_1,e_2)e_1=\alpha^2 e_2, \tilde{R}(e_1,e_3)e_1 = \tilde{R}(e_2,e_3)e_2 =\alpha^2 e_3, \tilde{R}(e_5,e_3)e_5=\alpha^2 e_3+\alpha e_4 $$
$$\tilde{R}(e_2,e_3)e_3= \tilde{R}(e_2,e_4)e_4=-\alpha^2 e_2, \tilde{R}(e_2,e_5)e_5=-\alpha^2 e_2+\alpha e_1 , \tilde{R}(e_3,e_4)e_4=-\alpha^2 e_3,$$
$$\tilde{R}((e_1,e_5)e_1= \tilde{R}(e_2,e_5)e_2= \tilde{R}(e_3,e_5)e_3= \tilde{R}(e_4,e_5)e_4=\alpha^2 e_5,$$
$$\tilde{R}(e_1,e_4)e_1= \tilde{R}(e_2,e_4)e_2= \tilde{R}(e_3,e_4)e_3=\alpha^2 e_4, \tilde{R}(e_5,e_4)e_5= \alpha^2 e_4-\alpha e_3.$$
The above expressions of $\tilde{R}$ satisfy our result \eqref{3.8} of Theorem (3.1).\\
Also from the above relations, we obtain,
\begin{equation}\label{5.1}
  \tilde{S}(e_1,e_1)= \tilde{S}(e_2,e_2)= \tilde{S}(e_3,e_3)= \tilde{S}(e_4,e_4)= \tilde{S}(e_5,e_5)=-4\alpha^2
\end{equation}
which satisfy our result \eqref{3.10} of Theorem (3.1).\\
Hence,
\begin{equation}\label{5.2}
  \tilde{r}=-20\alpha^2=r,
\end{equation}
which satisfies our result \eqref{3.11} of Theorem (3.1).\\
Now putting $X=Y=e_5$ in \eqref{4.5} and using \eqref{4.8}, we get,
\begin{equation}\label{5.3}
  (\tilde{\pounds}_\xi g)(e_5,e_5)=0.
\end{equation}
Then taking $X=Y=e_5$ in \eqref{4.3} and using \eqref{5.1}, \eqref{5.2}, \eqref{5.3}, we obtain,
\begin{equation}\label{5.4}
  -4\alpha^2\rho+\big[\Lambda+3\alpha^2\rho+2\alpha^2 q\big]+\big[\mu+\alpha^2\rho\big]=0,
\end{equation}
which gives,
\begin{equation}
  \Lambda+\mu=-2\alpha^2 q.\nonumber
\end{equation}
Hence the above equation proves that, $\Lambda$ and $\mu$ satisfies our result \eqref{4.9} of Theorem (4.1) for $n=5$ and so $g$ defines a $\ast$-$\eta$-Ricci-Yamabe soliton on the 5-dimensional $\alpha$-cosymplectic manifold with respect to a quarter-symmetric metric connection $\tilde{\nabla}$.

\section{Conclusion $\&$ Remarks}
The effect of $\ast$-$\eta$-Ricci-Yamabe soliton have been studied within the framework of $\alpha$-cosymplectic manifold. Here we have characterized $\alpha$-cosymplectic manifold, which admits $\ast$-$\eta$-Ricci-Yamabe soliton, in terms of quarter-symmetric metric connection. Also there is further scope of research in this direction of various types of complex manifolds like Kaehler manifolds, para-Kaehler manifolds, Hopf manifolds etc. There are some questions arise from our article to study further research.\\
(i) Are Theorem 4.1 true if we ignore quarter-symmetric connection on $\alpha$-cosymplectic manifold?\\
(ii) Are Theorem 4.4 true without assuming the potential vector field of the soliton is of gradient type?\\
(iii) If the conformal vector field is not killing, then is Theorem 6.6 true?\\
(iv) What are the results of the our paper also in true Trans-Sasakian manifolds, Co-k\"{a}hler manifold or para-contact geometry?

\end{document}